\def \be {\begin{equation}}
\def \bea {\begin{eqnarray}}
\def \ee {\end{equation}}
\def \eea {\end{eqnarray}}

\def \ll {\label}

\def \ox {\otimes}
\def \unit {{\bf 1}}
\def \complex {{\bf C}}
\def \newblock {}
\def \rf {(\ref}

\documentstyle[12pt,a4]{article}
\begin{document}
\bibliographystyle{unsrt}    % for BibTeX - sorted numerical labels
                 % by order of first citation.

\author{\LARGE \bf{\sc  Ladislav   Hlavat\'{y}}
\thanks{Postal   address:
B\v{r}ehov\'{a} 7, 115 19
Prague 1, Czech Republic. E-mail: hlavaty@br.fjfi.cvut.cz}
\\ {\it Department  of  Physics,}
\\ {\it  Faculty  of  Nuclear  Sciences  and
Physical Engineering}}
\title
{\LARGE \bf
{Solution of constant Yang--Baxter system in the dimension two}
}
\maketitle
q-alg/9511016
%\vskip 2cm
%Short title:{Solution of Yang-Baxter system}

%\vskip 1cm
%PACS numbers: 02.10 Tq, 03.65 Fd, 75.10 Jm.
%\newpage
\abstract{Complete solution of the
Yang--Baxter system (\ref{rrr0}) -- (\ref{zzr0}) in the dimension 2,
more precisely, all invertible $4\times 4$ matrices  $R,Q$ that
solve the system are given.}
\section{Introduction -- origin of the Yang--Baxter system}

The quantised braided groups were introduced recently in
\cite{hla:qbg}  %  (see also \cite{cou:qbg})
combining Majid's concept of braided groups \cite{maj:bg} and the FRT
formulation of quantum supergroups \cite{Liao}. The
 generators of quantised braided groups  $T=T_i^j, \  i,j \in
\{1,\ldots ,d=dimV\}$
satisfy the algebraic and braid relations
\begin{equation}
Q_{12}R_{12} ^{-1}T_{1}R_{12}T_{2}     =   R_{21}   ^{-1}
T_{2}R_{21}T_{1} Q_{12}             \label{qbgreln}
\end{equation}
\begin{equation}
\psi(T_{1}\ox R_{12}T_{2}  )   =   R_{12}T_{2}\ox R_{12}^{-1}
T_{1} R_{12}
\end{equation}
where the numerical matrices $Q,\ R$ satisfy the system of
Yang--Baxter--type equations
\begin{equation} Q_{12}Q_{13}Q_{23}     =
Q_{23}Q_{13}Q_{12},\label{rrr0}
\end{equation}
\begin{equation} R_{12}R_{13}R_{23}     =
R_{23}R_{13}R_{12},
\label{zzz0}
\end{equation}
\begin{equation} Q_{12}R_{13}R_{23}     =
R_{23}R_{13}Q_{12},\quad
\label{rzz0}
\end{equation}
\begin{equation}
 R_{12}R_{13}Q_{23}     =      Q_{23}R_{13}R_{12}.
\label{zzr0}
\end{equation}
The special cases of the quantised braided groups
are the quantum supergroups \cite{Liao}, quantum anyonic groups
\cite{majrod:afrt}
and braided $GL_q$ groups
\cite {cou:qbg}.

For obtaining all possible quantized braided groups
 we have to find all solutions of this Yang--Baxter system.
 (\ref{rrr0}) -- (\ref{zzr0}). Besides that solution of this
system give special solution of a more complicated system, see
e.g. \cite{vlad:qdouble}
\begin{equation} Q_{12}Q_{13}Q_{23}     =
Q_{23}Q_{13}Q_{12},\
\overline{Q}_{12}\overline{Q}_{13}\overline{Q}_{23}     =
\overline{Q}_{23}\overline{Q}_{13}\overline{Q}_{12},
\ll{vlad1}\end{equation}
\begin{equation} \overline{Q}_{12}R_{13}R_{23}     =
R_{23}R_{13} \overline{Q}_{12},\
 R_{12}R_{13}Q_{23}     =      Q_{23}R_{13}R_{12}.
\ll{vlad2}\end{equation}
Indeed, if $(Q,\, R)$ is a solution of
the system of Yang--Baxter type equations (\ref{rrr0}--\ref{zzr0})
then $(Q,\, R,\, \overline{Q}=PRPQR^{-1})$ is a solution of the
system \rf{vlad1},\ref{vlad2}). The latter system is, on the
other hand, a special case of  Yang--Baxter type equations
that appeared in \cite{frimai} as a consistency conditions for
quantization of nonultralocal models (see also
\cite{hla:afqnulm}).

There are several simple
solutions of the system (\ref{rrr0}--\ref{zzr0}).
One of them is $R=1,\ Q$ - any solution of the Yang--Baxter
equation (YBE). This gives the
algebras that correspond to the  ordinary  (unbraided)
quantum  groups in the FRT  formulation \cite{FRT}.
Rather trivial solution is $R,\ Q=P$ where P is the permutation
matrix and
$R$ any solution of the YBE. This solution does not yield  any
algebraic
relation because (\ref{qbgreln}) is satisfied identically in this
case.

Other simple solutions are  $(R,\,Q=R)$  or  $(R,\,Q=PR^{-1}P),\ R$
being any  solution  of  the  YBE. For $R$ that is both invertible and
has the so called second inversion $(R^{t_1})^{-1}$
they correspond    to    the (unquantised)    braided groups
\cite{maj:jmp91}.

In the paper \cite{hla:qbg} several nontrivial solutions of
(\ref{rrr0}) -- (\ref{zzr0}) were
presented but in general there is very little experience with
the Yang--Baxter systems to the best knowledge of the author.

{\em The goal of this paper  is} to give the complete solution of the
Yang--Baxter system (\ref{rrr0}) -- (\ref{zzr0}) in the dimension 2,
more precisely, {\em find all invertible $4\times 4$ matrices  $R,Q$
that
solve} (\ref{rrr0}) -- (\ref{zzr0}) .
The starting point for this is the Hietarinta's classification of
$4\times 4$ solutions of the YBE \cite{hie:ybecfnpl,hie:ybecfn}.

\section{Solution of the system}

Before  we start solving the system (\ref{rrr0}) -- (\ref{zzr0}) we
shall discuss its symmetries.
It is well known that each of the YBE (\ref{rrr0}) and (\ref{zzz0})
have
both continuous and discrete groups of symmetries. However, the
group of symmetries
of the system (\ref{rrr0}) -- (\ref{zzr0}) is not
the product  of the groups of symmetries  (\ref{rrr0}) and
(\ref{zzz0}).
On the other hand it is easy to check that
the system (\ref{rrr0}) -- (\ref{zzr0}) is invariant under
\be Q'= \lambda(S\ox S)Q(S\ox S)^{-1},
 \ R'=\kappa(S\ox S)R(S\ox S)^{-1},\ \lambda,\kappa\in \complex,\
S\in SL(2,\complex) \ll{sym1}\ee
and
\be Q''=Q^+=PQP,\  R''=R^+=PRP. \ll{sym2} \ee
They are the symmetries we are going to use.

The procedure we adopt for solving the system (\ref{rrr0}) --
(\ref{zzr0}) is the following.
{}From the Hietarinta's classification \cite{hie:ybecfnpl,hie:ybecfn}
we know that up to the symmetries
(\ref{sym1},\ref{sym2})
there are just eleven invertible, i.e. rank 4, solutions
of (\ref{zzz0}) (see also \cite{hla:usybe,hla:baxtn}). For each of
them
it is rather easy to solve the linear homogenous equations
(\ref{rzz0},\ref{zzr0}) for $Q$.
%We have used the computer algebra system
%REDUCE for this task.
Inserting the solution
of the linear equations into (\ref{rrr0}) we get equations for the
coefficients of linear combinations of matrices that solve
\rf{rzz0},\ref{zzr0}). Solving them we get the solution
of the whole system.

As mentioned in the introduction, it is clear that among others
we must obtain solutions $Q=P,\ Q=R,\ Q=PR^{-1}P$ (not necessarily
different) for each $R$. The interesting cases are those when some
other solutions appear.
As we shall see they are rather rare.

\subsection{Generic cases}
To demonstrate the method sketched above and show
a typical and simple example we choose
%the YBE solution %of \rf{zzz0}
\be R=\left( \begin{array}{cccc}
1&0&0&1\\0&1&1&0\\0&1&-1&0\\-1&0&0&1       \end{array}
\right)                              .\ee
This is the YBE
solution $R_{H0.2}$ in the   classification
\cite{hie:ybecfnpl,hie:ybecfn}
or $R_1$ in \cite{hla:usybe}. The solution of the linear equations
(\ref{rzz0},\ref{zzr0}) in this case is
\be Q=\alpha\left( \begin{array}{cccc}
1&0&0&0\\0&0&1&0\\0&1&0&0\\0&0&0&1      \end{array}
\right)+\beta\left( \begin{array}{cccc}
0&0&0&1\\0&1&0&0\\0&0&-1&0\\-1&0&0&0      \end{array}
\right)
\end{equation}
and from the YBE (\ref{rrr0}) we get cubic equation for the
coefficients $\alpha$
and $\beta$
\[ \beta(\beta^2-\alpha^2)=0. \]
It is easy to check that three solutions
$\beta=0,\pm\alpha$ of this equation correspond just to the expected
 solutions $Q=\alpha P,\ Q=\alpha R,\ Q=\alpha PR^{-1}P$ and no
other exists. These solution, where $R$ have the second inversion,
 yield (unquantized) braided group or
free algebras (for $Q=P$).

The same situation, namely that the only solutions of
\rf{rrr0})--\rf{zzr0}) are
$(R,\alpha P)$, $(R,\alpha R)$, $(R,\alpha PR^{-1}P)$, happens
{\em for all but  four invertible YBE solutions} $R$
of the Hietarinta's list.

The exceptional cases are investigated in the next section.
\subsection{Special cases}
The simplest solution of \rf{zzz0}) that admits
a wider set of solutions of the Yang--Baxter system
(\ref{rrr0},\ref{rzz0},\ref{zzr0})
than just $Q\propto P,R,PR^{-1}P$ is
\begin{equation}
R=\left( \begin{array}{cccc}
1&0&0&0\\0&1&0&0\\0&0&1&0\\1&0&0&-1       \end{array}
\right). \ \ee
This is a special case of the YBE
solution $R_{H1.2}$ in the   classification
\cite{hie:ybecfnpl,hie:ybecfn}
(or $R_{3}$ in \cite{hla:usybe}). %nontrivial solutions
 The solutions of the linear equations
(\ref{rzz0},\ref{zzr0}) in this case form four-dimensional space
\be Q=\left( \begin{array}{cccc}
\alpha+\beta&0&0&0\\0&\alpha-\delta&\beta+\delta&0\\0&\alpha&\beta&0\\
\gamma&0&0&\delta      \end{array}
\right)
\end{equation}
and from the YBE (\ref{rrr0}) we get three simple
cubic equations for the coefficients $\alpha$
, $\beta, \gamma$ and $\delta$
\be  \alpha\gamma(\beta+\delta)=0,\  \alpha\beta(\beta+\delta)=0,
\  \alpha(\alpha-\delta)(\beta+\delta)=0.\ee
The invertible solutions of these equations up to the

transformation are $Q=\alpha P$ or
\be Q=\beta\left( \begin{array}{cccc}
1&0&0&0\\0&1&0&0\\0&1-t&t&0\\r&0&0&-t       \end{array}
\right). \ \ll{r3}\ee
Note that in general, the matrix $Q$ %given by (\ref{r3})
is proportional
neither to $R$ nor $PR^{-1}P$.

The second type of nontrivial solutions to the system
(\ref{rrr0}) -- (\ref{zzr0}) is obtained when
\be R=\left( \begin{array}{cccc}
0&0&0&1\\0&0&t&0\\0&t&0&0\\1&0&0&0       \end{array}
\right)                              .\ \ee
This is the
solution $R_{H1.4}$ in the Hieterinta's classification
\cite{hie:ybecfnpl,hie:ybecfn}
or $R_9$ in \cite{hla:usybe}. The solution space of the linear
equations
(\ref{rzz0},\ref{zzr0}) in this case is three-dimensional
\be Q=\left( \begin{array}{cccc}
\alpha&0&0&\gamma\\0&0&\beta&0\\0&\beta&0&0\\ \gamma&0&0&\alpha
\end{array}
\right)
\end{equation}
and two cubic equations for the coefficients $\alpha$
, $\beta$ and $\gamma$ are obtained from the YBE (\ref{rrr0})
\be \alpha^2 \gamma=0,\ \  \alpha(\gamma^2+\alpha\beta-\beta^2)=0. \ee
They yield $Q=\alpha P$,
\be Q=\gamma\left( \begin{array}{cccc}
0&0&0&1\\0&0&t'&0\\0&t'&0&0\\1&0&0&0       \end{array}
\right)                              \ll{r9} \ee
or nonivertible matrices. Once again
 the matrix $Q$ given by (\ref{r9}) is proportional
neither to $R$ nor $PR^{-1}P$, in general.

The third class of nontrivial
%YBE solution that provides a wider set of
solutions of the Yang--Baxter system (\ref{rrr0}) -- (\ref{zzr0})
come from
\begin{equation}
R=\left( \begin{array}{cccc}
1&0&0&0\\x&1&0&0\\y&0&1&0\\z&y&x&1       \end{array}
\right). \ \ee
This is the YBE
solution $R_{H2.3}$ in the   classification
\cite{hie:ybecfnpl,hie:ybecfn}
or $R_{10}$ in \cite{hla:baxtn}. The solutions of the linear equations
(\ref{rzz0},\ref{zzr0}) in this case form six-dimensional space
\be Q=\left( \begin{array}{cccc}
\alpha_1&0&0&0\\ \beta_1&\alpha_2&\alpha_1-\alpha_2&0\\
\beta_2&\alpha_1-\alpha_2&\alpha_2&0\\ \gamma&\delta_1&
\beta_1+\beta_2-\delta_1&\alpha_1
 \end{array}
\right) \ \ee
%\[\beta_1+\beta^2=\delta_1+\delta_2 \]
and  solutions of
the set of cubic equations for the coefficients $\alpha,\ldots,
\delta_1$
that we get from the YBE (\ref{rrr0}) give two nonivertible matrices
$Q$.
\be Q=\left( \begin{array}{cccc}
1&0&0&0\\a&1&0&0\\b&0&1&0\\c&b&a&1       \end{array}
\right) \ \ee
and
\be Q=\left( \begin{array}{cccc}
1&0&0&0\\-g&1&0&0\\g&0&1&0\\-gh&h&-h&1       \end{array}
\right).
\end{equation}
These solutions were found in \cite{hla:qbg} and the corresponding
quantized braided groups were presented.

The last special case
is given by diagonal $R$. This possibility
 was investigated in the paper \cite{hla:qbg}
and the following lemma was proved.\\
%so that we cite only the results here:              \\
{\bf Lemma}: If $R$ is a diagonal matrix, $R =
diag(x,u,v,y),\ xuvy  \neq  0  $,
then    there    are    three    types    of     solutions     of
the system (\ref{rrr0}) -- (\ref{zzr0}) :
\\
1) $R\propto\unit$ i.e. $x=u =v=y$ and $ Q$ is an arbitrary YBE
solution.
\\
2) $R$ is proportional to sign--diagonal, i.e. $x^{2}= u^{2}= v^{2}=
y^{2}$,
and Q is an arbitrary  YBE
 solution of the eight--or--less--vertex form
\begin{equation}
Q=
\left( \begin{array}{cccc}
q&0&0&a\\0&r&b&0\\0&c&s&0\\d&0&0&t       \end{array}
\right)
\label{8vf} \end{equation}
\\
3)$R$ is general diagonal matrix and $Q$ is an arbitrary six--or--
less--vertex
YBE solution  i.e. of  the  form
(\ref{8vf})  where
$a=d=0$.

When classifying solutions following from these three
possibilities we must remember  that in general
we do not have at our disposal the symmetry
transformation
\be Q'=(S\ox S)Q(S\ox S)^{-1},\ S\in SL(2,\complex)
 \ll{symr} \ee
used in the classification of the YBE solutions
because now it must be accompanied by
\be  \ R'=(S\ox S)R(S\ox S)^{-1} \ll{symz} \ee
that might have been used for bringing $R$
to the diagonal form. It means
that we can use only such transformation of $Q$ that keep $R$
diagonal.

For scalar $R$ i.e. $R=\lambda \unit$ this is no restriction
so that we get as many nonequivalent
solutions of the system (\ref{rrr0}) -- (\ref{zzr0}) as there
are solutions of the YBE \cite{hie:ybecfnpl,hie:ybecfn}.

For the other types
of diagonal $R$,
it rather easy to see that %such transformation are
 we can use only transformations
\rf{symr},\ref{symz}) %of $Q$ that
generated either by diagonal or antidiagonal matrix $S$. It means
that we must
classify the corresponding type of the YBE solutions up to this
special type of symmetries.

Classification of the six--or--less vertex solutions up to this
restricted type of transformations is not different from
the usual one \cite{hla:usybe}. There are just four
invertible solutions of the YBE, namely
\begin{equation}
\left( \begin{array}{cccc}
q&0&0&0\\0&1&0&0\\0&q-t&qt&0\\0&0&0&q      \end{array}
\right)                            ,\
\ \left( \begin{array}{cccc}
q&0&0&0\\0&1&0&0\\0&q-t&qt&0\\0&0&0&-t      \end{array}
\right)                            ,\
\ \left( \begin{array}{cccc}
a&0&0&0\\0&b&0&0\\0&0&c&0\\0&0&0&d       \end{array}
\right)                              ,
\end{equation}
and the permutation matrix.%$$Q=xP$.
%The  numbering  corresponds  to  that in classifications   given   in
%\cite{hla:usybe}.

On the other hand, classification of the invertible eight--or--less
vertex
YBE solutions \rf{8vf})
up to the restricted transformations is different
from the usual one because besides the
solutions
%(\ref{z1}, \ref{r3}, \ref{r9}, \ref{6v}) and
presented in the classification \cite{hla:usybe}, there are YBE
solutions
\be Q=\left( \begin{array}{cccc}
a&0&0&x\\0&\pm a&x&0\\0&x&\pm a&0\\x&0&0&a      \end{array}
\right),\ \ x,a\neq 0                             \ll{spec8v}
\ee
that can be transformed %by (\ref{}) with (anti)diagonal $S$ to
solution
%in the usual classification
to those in the paper \cite{hla:usybe} only by non-(anti)diagonal
transformations
(\ref{symr},\ref{symz}).

%To summarize the fourth special case, we can say that pairs of
%invertible matrices ($R,Q$) where $R$ is diagonal
%solve the system (\ref{rrr0}) -- (\ref{zzr0}) if and only if
%$Q$ is
%one of \rf{6v}) or $Q=xP$ or $R$ is proportional to the
%sign--diagonal matrix, i.e. \rf{sd}) holds and $Q$ is one of the
%matrices  (\ref{r1}, \ref{r3}, \ref{r9}, \ref{6v}, \ref{Q27}),
%\ref{spec8v}) or $R=x\unit$ and

\section{Conclusions}

We have classified all invertible solutions
of the Yang--Baxter system (\ref{rrr0}) -- (\ref{zzr0})
up to the transformations \rf{sym1},\ref{sym2}).
%with $Q$ and $R$ invertible.
They are given by pairs ($R,Q$)
where $R$ belong to the Hietarinta's list
\cite{hie:ybecfnpl,hie:ybecfn}.
%We have found that generically the ordinary solutions i.e.

For all but four types of invertible solutions from the Hietarinta's
list, the only solutions of the Yang--Baxter system (\ref{rrr0}) --
(\ref{zzr0})
are pairs
$(R,\lambda R)$ or
$(R,\lambda PR^{-1}P)$ or
$(R,\lambda P)$
 where %$R$ is a solution of the YBE and
$P$ is the permutation matrix.

The list of exceptional solutions of the Yang--Baxter system
(\ref{rrr0}) -- (\ref{zzr0})
 up to the transformations \rf{sym1},\ref{sym2}) is
\begin{equation}
R=\left( \begin{array}{cccc}
1&0&0&0\\0&1&0&0\\0&0&1&0\\1&0&0&-1       \end{array}
\right),\ \
 Q=\left( \begin{array}{cccc}
1&0&0&0\\0&1&0&0\\0&1-t&t&0\\r&0&0&-t       \end{array}
\right),\ee
\be R=\left( \begin{array}{cccc}
0&0&0&1\\0&0&t&0\\0&t&0&0\\1&0&0&0       \end{array}
\right),\
Q=\left( \begin{array}{cccc}
0&0&0&1\\0&0&a&0\\0&a&0&0\\1&0&0&0       \end{array}
\right) ,                             \ll{cr9} \ee
\begin{equation}
R=\left( \begin{array}{cccc}
1&0&0&0\\x&1&0&0\\y&0&1&0\\z&y&x&1       \end{array}
\right), \  Q=\left( \begin{array}{cccc}
1&0&0&0\\a&1&0&0\\b&0&1&0\\c&b&a&1       \end{array}
\right), \ \ee
\be
R=\left( \begin{array}{cccc}
1&0&0&0\\x&1&0&0\\y&0&1&0\\z&y&x&1       \end{array}
\right), \ Q=\left( \begin{array}{cccc}
1&0&0&0\\-g&1&0&0\\g&0&1&0\\-gh&h&-h&1       \end{array}
\right).
\end{equation}
\be R=\left( \begin{array}{cccc}
1&0&0&0\\0&x&0&0\\0&0&y&0\\0&0&0&z       \end{array}
\right),\ \ Q=\left( \begin{array}{cccc}
q&0&0&0\\0&1&0&0\\0&q-t&qt&0\\0&0&0&q      \end{array}
\right)                            ,\ee
\begin{equation}
R=\left( \begin{array}{cccc}
1&0&0&0\\0&x&0&0\\0&0&y&0\\0&0&0&z       \end{array}
\right),\ \ Q=\left( \begin{array}{cccc}
q&0&0&0\\0&1&0&0\\0&q-t&qt&0\\0&0&0&-t      \end{array}
\right)                            ,\ee
\begin{equation}
R=\left( \begin{array}{cccc}
1&0&0&0\\0&x&0&0\\0&0&y&0\\0&0&0&z       \end{array}
\right),\ \ Q=\left( \begin{array}{cccc}
1&0&0&0\\0&a&0&0\\0&0&b&0\\0&0&0&c       \end{array}
\right)                              ,\ll{6v}
\end{equation}
and pairs $(R,Q)$ where
\be R=\left( \begin{array}{cccc}
1&0&0&0\\0&x&0&0\\0&0&y&0\\0&0&0&z       \end{array}
\right),\ x^2=y^2=z^2=1\ee
and $Q$ is any solution of the following list
\be Q=\left( \begin{array}{cccc}
1&0&0&1\\0&1&1&0\\0&1&-1&0\\-1&0&0&1       \end{array}
\right)
,\ \  Q=\left( \begin{array}{cccc}
1+t&0&0&1\\0&s&1&0\\0&1&s&0\\1&0&0&1-t      \end{array}
\right), \ s^2=1+t^2,
\ee
\be Q=\left( \begin{array}{cccc}
1&0&0&0\\0&1&0&0\\0&1-t&t&0\\1&0&0&-t       \end{array}
\right),\ \
Q=\left( \begin{array}{cccc}
1&0&0&0\\0&-1&0&0\\0&0&-1&0\\1&0&0&1       \end{array}
\right)
\ee
\be Q=\left( \begin{array}{cccc}
0&0&0&1\\0&0&t&0\\0&t&0&0\\1&0&0&0       \end{array}
\right)
,\ \ Q=\left( \begin{array}{cccc}
a&0&0&1\\0&b&1&0\\0&1&b&0\\1&0&0&a      \end{array}
\right),\ b^2=a^2.
\ee
This list is exhaustive i.e. any other invertible solution of the
system
  (\ref{rrr0}) -- (\ref{zzr0})
can be obtained by the transformations \rf{sym1},\ref{sym2}).

The author gratefully acknowledges the support of the grant
No. 202/93/1314 of the Czech Republic and the grant No.8154 of
the Czech Technical University.

\end{document}